# Polyhedral Tori with Minimal Coordinates

## Authors

Stefan Hougardy, Frank H. Lutz, and Mariano Zelke

## Description

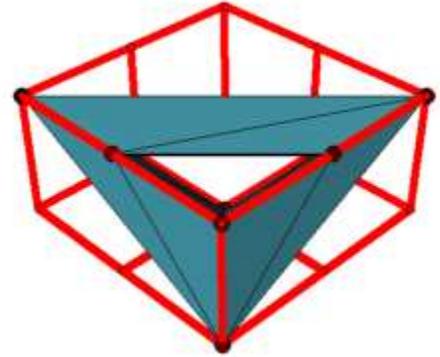

We give explicit realizations with small integer coordinates for all triangulated tori with up to 12 vertices. In particular, we provide coordinate-minimal realizations in general position for all triangulations of the torus with 7, 8, 9, and 10 vertices. For the unique 7-vertex triangulation of the torus we show that all corresponding 72 oriented matroids are realizable in the 6x6x6-cube. Moreover, we present polyhedral tori with 8 vertices in the 2x2x2-cube, general position realizations of triangulated tori with 8 vertices in the 2x2x3-cuboid as well as polyhedral tori with 9 and 10 vertices in the 1x2x2-cuboid.

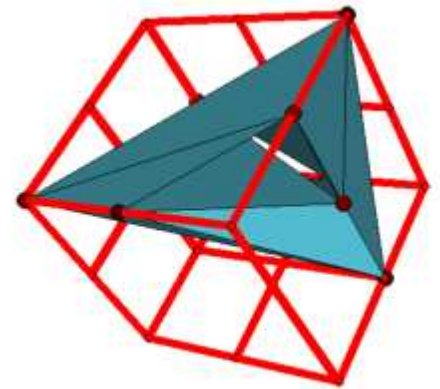

According to Heawood's bound [12], every triangulation of a (closed) surface M of Euler characteristic chi(M) has at least

    n ≥ 1/2(7+sqrt(49-24*chi(M)))

vertices. Thus, at least 7 vertices are necessary to triangulate the 2-dimensional torus $\mathbf{T}^2$ with chi($\mathbf{T}^2$)=0. The unique vertex-minimal triangulation with 7 vertices of the torus was first described by Möbius [20], a first geometric realization of it in $\mathbf{R}^3$ was given by Császár [8]; see also [10, 16].

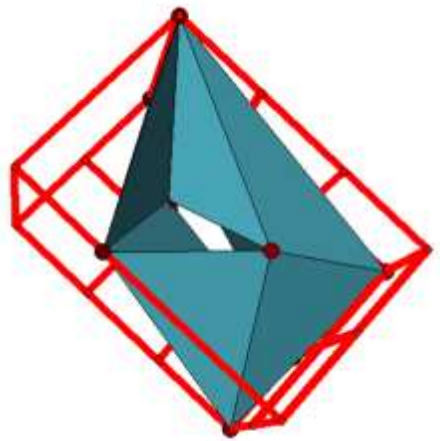

By Steinitz' theorem (cf. [24, Ch. 4]), every triangulated 2-sphere is realizable as the boundary complex of a convex 3-dimensional polytope. For orientable surfaces of genus g ≥ 1 Grünbaum [11, Ch. 13.2] asked whether triangulations of these surfaces can always be *realized geometrically as polyhedra* in $\mathbf{R}^3$, i.e., with straight edges, flat triangles, and without self intersections? Recently, Grünbaum's question was settled for the torus by Archdeacon, Bonnington, and Ellis-Monaghan [1]: All triangulated tori *are* realizable. In general, however, the answer to Grünbaum's problem turned out to be "*No*": Bokowski and Guedes de Oliveira [5] showed that there is a non-realizable triangulation of the orientable surface of genus g = 6. Schewe [21] extended this result by proving that there are non-realizable triangulations for all orientable surfaces of genus g ≥ 5. For surfaces of genus 2 ≤ g ≤ 4 the problem remains open. However, all vertex-minimal triangulations of these surfaces are realizable [3, 13, 14, 15, 17].

The realizability problem for triangulated surfaces is decidable (cf. Bokowski [2] and Bokowski and Sturmfels [7, Ch. VIII]), but there is *no* deterministic algorithm known that would solve the realization problem for instances with as few as, say, 7 vertices in reasonable time. At present, however, there are three heuristic procedures available to search for realizations:

(1) by minimizing the intersection edge functional [15],

(2) by Bokowski's rubber band method [3], and

(3) by searching for small coordinates [13, 14].

In fact, there are three basic types of realizations:

*(i) linear embeddings* or *realizations* with straight edges, flat triangles, and without self intersections,

*(ii) proper realizations* which are realizations with no coplanar neighboring faces, and

*(iii) realizations in general position* with no three vertices on a line and no four vertices on a plane.

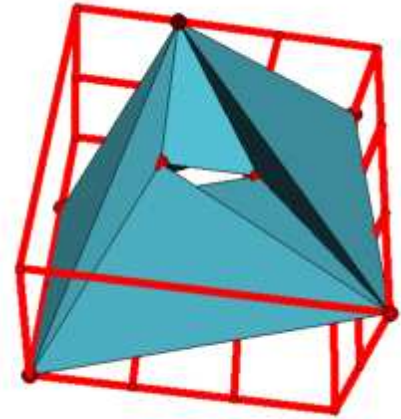

The combinatorics of a realization is encoded by an *oriented matroid*; cf. [7]. For general position (gp) realizations the corresponding oriented matroids are *uniform*. If a triangulation is *neighborly*, i.e., if the 1-skeleton of the triangulation is a complete graph, then every realization of it is in general position (since otherwise some two edges would intersect).

The unique triangulation of the torus with 7 vertices is neighborly, and there are (up to equivalence) 72 different oriented matroids that are compatible with this 7-vertex triangulation [4]. (The listing of these 72 oriented matroids in [4] contains two typos: the signs at the third to last positions of b16 and b18 have to be switched.) It is known that all (rank four) oriented matroids with up to 8 elements are realizable [6], and therefore all the 72 Möbius torus oriented matroids are realizable. We used technique (3) from [13, 14] to obtain realizations for the 72 oriented matroids with small integer coordinates.

**Theorem:** Of the 72 Möbius torus oriented matroids 13 have coordinate-minimal integer realizations in the 2x3x3-cuboid, 18 in the 3x3x3-cube, 32 in the 4x4x4-cube, and 7 in the 5x5x5-cube, respectively. The two remaining of the 72 oriented matroids are realizable in the 6x6x6-cube, but are not realizable in the 4x4x4-cube.

Altogether, there are (up to symmetry) 46 combinatorially distinct coordinate minimal realizations of the Möbius torus in the 2x3x3-cuboid. (Differing from [4], we used the lexicographic minimal vertex labeling for the Möbius torus. A mapping between the vertex labeling from [4] and the lexicographic minimal labeling is given by the permutation (1,2,5,4,6,7).) One of the realizations of the Möbius torus in the 2x3x3-cuboid is displayed in the last applet.

A complete enumeration of all triangulated surfaces with up to 12 vertices was obtained in [17] and [23]; see [18] for a list of facets (in mixed lexicographic format as well as in lexicographic format). In particular, there are 7 triangulations of the torus with 8 vertices, 112 with 9, 2109 with 10, 37867 with 11, and 605496 with 12 vertices, respectively. With his program surftri [22], Sulanke even obtained all triangulated tori with up to 17 vertices.

**Theorem:** There are no general position realizations of triangulated tori in the 2x2x2-cube.

*Proof:* The 2x2x2-cube consist of three layers of size 2x2. By the general position requirement, each of the three layers can contain at most three vertices. However, by our enumeration, there is no general position realization in the 2x2x2-cube of a triangulated torus with at most 9 vertices.

**Theorem:** There are exactly two triangulated tori, with 8 vertices, that have general position realizations in the 2x2x3-cuboid.

The third applet displays one of the two general position realizations in the 2x2x3-cuboid.

In a non-proper linear realization, two neighboring coplanar triangles either form a square or a larger triangle (see, for example, the first or the second applet). If all maximal unions of neighboring coplanar triangles form convex polygons and if all resulting convex polygons intersect in edges and vertices only, then the linear realization of the triangulation induces a proper realization of a polyhedral map with fewer edges and possibly fewer vertices and non-triangular faces.

According to our search, there are no general position realizations of the Möbius torus in the 2x2x2-cube. However, we were able to find linear realizations with coplanar neighboring triangles for 5 of the 7 triangulations of the torus with 8 vertices in the 2x2x2-cube. Moreover, 31 of the 112 triangulations of the torus with 9 vertices as well as 567 of the 2109 triangulations of the torus with 10 vertices even have linear realizations in the 1x2x2-cuboid. The first applet presents one of the examples with 9 vertices in the 1x2x2-cuboid, the second applet presents one of the examples with 8 vertices in the 2x2x2-cube.

**Theorem:** There are linear realizations of triangulated tori in the 1x2x2-cuboid.

These realizations are smallest possible.

**Theorem:** There are no linear realizations of triangulated tori in the 1x1x2-cuboid.

*Proof:* The 1x1x2-cuboid contains 12 integer vectors. However, not all of the 66 edges connecting these points can be present at the same time in a corresponding linear realization. For every square of size 1x1 on the boundary, only one of the diagonals can be present, which forbids 10 edges. In addition, in the quadrilaterals on the boundary of size 1x2 at most one of the diagonals or the middle edge can be present, ruling out 8 further edges. Each of the two 1x1x1-cubes has four space diagonals of which only one can be present, respectively, ruling out 6 edges. The center point of the 1x1x2-cuboid lies on 6 edges, which forbids 5 of them. Finally, the edges of length 2 on the boundary can be split into two edges of length 1, but are in conflict with these, thus reducing the admissible number of edges by additional 4 edges. In conclusion, at most 33 of the 66 edges can be present at the same time. A triangulated torus with n vertices has 3n edges. Thus, a realization, if one exists, cannot use all the 12 vertices. But if one of the vertices is not used, then less than the 33 of the 66 edges can be present at the same time. Therefore, there are no polyhedral tori in the 1x1x2-cuboid with 11 vertices. By enumeration, the triangulations of the torus with 8, 9, and 10 vertices do not have linear realizations in the 1x1x2-cuboid. The Möbius torus with 7 vertices is neighborly and therefore also cannot be realized in the 1x1x2-cuboid.

**Theorem:** All 2229 triangulations of the torus with up to 10 vertices are realizable in general position in the 4x4x4-cube, with 2218 of the examples even having general position realizations in the 3x3x3-cube. The remaining 11 examples, all with 10 vertices, cannot be realized in general position in the 3x3x3-cube.

For the triangulated tori with up to 10 vertices we used technique (3) from [13, 14] to obtain the above general position realizations with small respectively minimal integer coordinates. Fendrich [9] provided realizations for all triangulated tori with up to 11 vertices (most of them via embeddings in the 2-skeleta of random 4-polytops, cf. [3]). By using technique (1) from [15] we obtained small integer realizations in general position for the triangulated tori with 11 and 12 vertices. It took 345 CPU minutes (on a 3.5 GHz prozessor) to find realizations for the 37867 triangulated tori with 11 vertices in the 20x20x20-cube. Hereby, we used *recycling of coordinates* (cf. [15, 17]): In the lexicographically sorted list of triangulations two consecutive triangulations often differ only slightly. Whenever we find coordinates for one example from the list we reuse these coordinates as starting coordinates for the next example from the list. In many cases, only few coordinate moves then are necessary to obtain a realization. An independent run of the 37867 triangulated tori with 11 vertices in the 40x40x40-cube was slightly faster, this time taking 317 CPU minutes. For the 605496 tori with 12 vertices about 70 CPU hours were necessary in total to obtain realizations in the 40x40x40-cube.

All the realizations of triangulated tori with up to 12 vertices are unknotted. In fact, at least 16 vertices are necessary for a knotted polyhedral torus; see [19].

**Index of applets**

Applet 1: Linear embedding of a triangulated torus with 9 vertices in the 1x2x2-cuboid.

Applet 2: Linear embedding of a triangulated torus with 8 vertices in the 2x2x2-cube.

Applet 3: General position realization of a triangulated torus with 8 vertices in the 2x2x3-cuboid.

Applet 4: General position realization of the Möbius torus with 7 vertices in the 2x3x3-cuboid.

## Acknowledgement

The first and the third author were supported by the DFG Research Center MATHEON "Mathematics for Key Technologies", Berlin, the second author was supported by the DFG Research Group "Polyhedral Surfaces", Berlin.

**Keywords**               triangulated surface; polyhedral realization; small coordinates

**MSC-2000 Classification**   52B70 (57Q15)

**Authors' Addresses**

Stefan Hougardy

> Universität Bonn
> Forschungsinstitut für Diskrete Mathematik
> Lennéstr. 2
> 53113 Bonn
> Germany
> hougardy@or.uni-bonn.de
> http://www.or.uni-bonn.de/cards/card-hougardy.de.html

Frank H. Lutz

> Technische Universität Berlin
> Institut für Mathematik, Sekr. MA 3-2
> Straße des 17. Juni 136
> 10623 Berlin
> Germany
> lutz@math.tu-berlin.de
> http://www.math.tu-berlin.de/~lutz

Mariano Zelke

> Humboldt-Universität zu Berlin
> Institut für Informatik
> Unter den Linden 6
> 10099 Berlin
> Germany
> zelke@informatik.hu-berlin.de
> http://www.informatik.hu-berlin.de/~zelke/